\documentclass[a4paper,12pt]{amsproc}
\usepackage{amsmath,amsfonts,amssymb,amsthm,anysize}%,dsfont}
\usepackage{enumerate}
\usepackage{epsfig}
\usepackage{supertabular}
\usepackage{graphicx}

\newcommand\Z{{\mathbb Z}}

\newcommand\F{{\mathbb F}}

\newcommand\AG{\mathsf{AG}}

\renewcommand\mod{{\mathrm{mod\, \, }}}

\theoremstyle{plain}
\newtheorem{theorem}{Theorem}[section]

\newtheorem{lemma}[theorem]{Lemma}

\newtheorem{proposition}[theorem]{Proposition}

\numberwithin{equation}{section}
\theoremstyle{remark}

\usepackage{url, hyperref}
\hypersetup{citecolor=blue, linkcolor=blue, colorlinks=true}

\begin{document}

\title[A new family of Hadamard matrices of order $4(2q^2+1)$]{A new family of Hadamard matrices of order $4(2q^2+1)$}

\author{Ka Hin Leung}
\address{ %
Department of Mathematics\\
National University of Singapore, 
Kent Ridge, Singapore 119260, Republic of Singapore}
\email{matlkh@nus.edu.sg}

\thanks{}

\author{Koji Momihara}
\address{ %
Division of Natural Science,\\
Faculty of Advanced Science and Technology,\\
Kumamoto University\\
2-40-1 Kurokami, Kumamoto 860-8555, Japan}
\email{momihara@educ.kumamoto-u.ac.jp}
\thanks{Koji Momihara was supported by 
JSPS under Grant-in-Aid for Young Scientists (B) 17K14236 and Scientific Research (B) 15H03636.}

\author{Qing Xiang}
\address{ %
Department of Mathematical Sciences\\
University of Delaware\\
Newark DE 19716, USA
}
\email{qxiang@udel.edu}
\thanks{Qing Xiang was supported by an NSF grant DMS-1600850.}

 %(primary), 05C50, 05B10, 11T24 (secondary)
%}

\begin{abstract}
Let $q$ be a prime power of the form $q=12c^2+4c+3$ with $c$ an arbitrary integer. In this paper we construct a difference family with parameters $(2q^2;q^2,q^2,q^2,q^2-1;2q^2-2)$ in $\Z_2\times (\F_{q^2},+)$. As a consequence, by applying the Wallis-Whiteman array, we obtain Hadamard matrices of order $4(2q^2+1)$ for the aforementioned $q$'s.   
\end{abstract}

% 2000 MSC numbers
% 05B25 Finite geometries

% 05E30 Association schemes, strongly regular graphs
% 51E12 Generalized quadrangles, generalized polygons

\maketitle

%%%%%%%%%%%%%%%%%%%%%%%%%%%%%%%%%%%%%%%%%%%%%%%%%%%%%%%%%%%
%%%%%%%%%%%%%%%%%%%%%%%%%%%%%%%%%%%%%%%%%%%%%%%%%%%%%%%%%%%
\section{Introduction}\label{sec:HP}
%This is a continuous paper of \cite{LM18}. 
A {\it Hadamard matrix} of order $v$ is a $v\times v$ matrix $H$ with entries $\pm 1$ such that $HH^{\top} = vI$, where $I$ is the identity matrix. It can be easily shown that if $H$ is a Hadamard matrix of order $v$, then $v=1, 2$, or $4t$ for some positive integer $t$. A long-standing conjecture in combinatorics states that a Hadamard matrix of order $v$ exists for every $v\equiv 0\,(\mod{4})$.  Despite the work of many researchers, the conjecture is far from being resolved. Currently it is still not known whether the set of orders of Hadamard matrices has positive density. For some sparse infinite subsequences of $\{4t: t=1,2,3,\ldots\}$, it is often possible to construct Hadamard matrices of order $v$ for every $v$ belonging to the subsequences. The most famous examples are the Paley constructions which produce Hadamard matrices of order $q+1$ if $q$ is a prime power congruent to 3 modulo 4, and Hadamard matrices of order $2(q+1)$ if $q$ is a prime power congruent to 1 modulo 4. As further examples, we mention that for prime powers $q\equiv 1\,(\mod{4})$ or $q\equiv 3\,(\mod{8})$, Xia and Liu \cite{XL91, XL96} construct Hadamard matrices of order $4q^2$; for $q\equiv 7\,(\mod{8})$, the first author, Ma and Schmidt \cite{LMS06} construct two possibly infinite families of Hadamard matrices of order $4q^2$. All these constructions are based on cyclotomy of finite fields. The Paley constructions use the nonzero squares of $\F_q$. The constructions by Xia and Liu \cite{XL91, XL96}, and by Leung, Ma, and Schmidt \cite{LMS06} use the $4^{\rm th}$, $8^{\rm th}$ and $(q+1)^{\rm th}$ cyclotomic classes of  $\F_{q^2}$. The main idea behind the constructions of Xia/Liu and Leung/Ma/Schmidt is to use cyclotomic classes of finite fields to construct a difference family with appropriate parameters in an abelian group $G$.

Throughout this paper,  we will use the following notation. Let $(G,+)$ be an additively written finite abelian group and  let $G^\ast:=G\setminus \{0_G\}$. For any subset $D$ in $G$, we define 
$D^{(-1)}:=\{-x:x \in D\}$, 
$\overline{D}:=G^\ast \setminus D$, and $D^c:=G \setminus D$. 
Furthermore, we will identify $D$ with the group ring element $\sum_{x\in D}x\in \Z[G]$ when there is no confusion.

Let $B_i$, $i=1,2,\ldots,\ell$, be $k_i$-subsets of $G$.
%Any subset $A$ of $G$ is said to be {\it symmetric} if $ A=A^{(-1)}$.  A family ${\mathcal B}$ is said to be {\it symmetric} if each $B_i$ is symmetric.  
The set ${\mathcal B}=\{B_i:i=1,2,\ldots,\ell\}$ is called a 
{\it difference family with parameters $(v;k_1,k_2,\ldots,k_\ell;\lambda)$ in $G$} if 
the list of differences ``$x-y, x,y\in B_i,x\not=y,i=1,2,\ldots,\ell$" represents every nonzero element of $G$ exactly $\lambda$ times; or equivalently  
\[ \sum_{i=1}^{\ell} B_i B_i^{(-1)}= \lambda G+ \Big(\sum_{i=1}^{\ell} k_i-\lambda\Big) \cdot 0_G.\]
%\begin{remark}
%To avoid possible confusion, we write $\lambda \cdot 0_G$ as the element in the group ring such that the coefficient of $0_G$ is $\lambda$. This notation will be used throughout this paper. 
%\end{remark}
Each subset $B_i$ is called 
a {\it block} of ${\mathcal B}$. %If there is only one block in a difference family, the block is a {\it difference set}. 
We now define two special classes of difference families.  
A difference family in $G$ with four blocks is said to be of {\it type $H$} if $\sum_{i=1}^4 k_i-|G|=\lambda$; and  of 
{\it type $H_4^\ast$} %difference family with $\ell$ blocks, 
if  $\sum_{i=1}^4 k_i-(|G|+1)=\lambda$. 

%Let us first summarize how a difference family be used to construct a %Hadamard matrix. 
It is well known that if there is a difference family of type $H$ in $G$, then we obtain a Hadamard matrix of order $4|G|$ by plugging the group invariant  %circulant ({\color{red}why circulant here? The group $G$ is only abelian, it may not be cyclic. Should it be group invariant instead of circulant?}) 
$(-1,1)$ matrices obtained from its blocks into %the Williamson array~\cite{W44,W47} or 
the Goethals-Seidel array~\cite{GS67}. In the literature, difference families of type $H$ have been extensively studied~\cite{LMS06,XL91,XL95,XL96,XX94,XX99,XXSW05}.

On the other hand, from a difference family of type $H_4^\ast$ in a finite abelian group $G$,  we obtain a Hadamard matrix of order $4(|G|+1)$ by plugging the the group invariant $(-1,1)$ matrices obtained from its blocks into %the Williamson array~\cite{W44,W47} or 
the Wallis-Whiteman array~\cite[Theorem~4.17]{WSW}. %In the case where $\ell=8$, we may use the Kharaghani array under the assumption that the circulant $(-1,1)$ matrices $M_i$, $i=0,1,\ldots,7$, obtained from its blocks are amicable, i.e., $\sum_{i=0}^3
%(M_{2i}M_{2i+1}^\tra-M_{2i+1}M_{2i}^\tra)=O$. (See, e.g., 
%\cite{K00},\cite[Lemma~4.20]{S17},\cite[Page~12]{SY89}.)  Note that the amicability condition becomes trivial if the difference family is symmetric. That explains why we are interested in difference families of type $H_4^\ast$ and symmetric difference families of type $H_8^\ast$ in this paper. 
Indeed difference families of type $H_4^\ast$ are particularly interesting as the orders of the Hadamard matrices obtained from the difference families are no longer of the form $4|G|$, but of the form $4(|G|+1)$. Very recently, the first and second authors \cite{LM18} gave two new constructions of difference families of type $H_4^\ast$ with parameters $(2n;n,n,n,n-1;2n-2)$. Difference families with these parameters were initially considered by Whiteman~\cite{W72}, who obtained one infinite family. Soon afterwards,
Spence~\cite{S75} came up with two new families whose constructions are based on relative difference sets.  On the other hand, the existence of difference families with parameters $(2n;n,n,n,n-1;2n-2)$ in dihedral groups was also studied in \cite{K96,KN02,SY00}. Let us summarize all known constructions of difference families of type $H_4^\ast$ with parameters  $(2n;n,n,n,n-1;2n-2)$.

\begin{theorem}
%Let $G$ be an abelian group of order $n$.  
There exists a difference family of type $H_4^\ast$ with parameters $(2n;n,n,n,$ $n-1;2n-2)$ %in $\Z_{2}\times G$ 
if $n$ satisfies any of the following conditions: 
\begin{itemize}
\item[(1)] {\em \cite{W72,SY00}} $n=q$ and $2q-1$ are both prime powers.
\item[(2)] {\em \cite{S75}} $q=2n+1$ is a prime power for which there exists a nonnegative integer $s$ such that $(q-2^{s+1}-1)/2^{s+1}$ is an odd prime power.
\item[(3)] {\em \cite{S75}} $n=q$ is a prime power such that $q\equiv 1\,(\mod{4})$, and $q-2$ is also a  prime power. 
\item[(4)] {\em \cite{LM18}} $n=9^{t_0}q_1^{4t_1}q_2^{4t_2}\cdots q_s^{4t_s}$, where 
$p_i$, $i=1,2,\ldots,s$, are prime powers and $t_i$, $i=0,1,\ldots,s$, are  nonnegative integer.
\item[(5)] {\em \cite{LM18}} $n=q^2$ with $q$ a prime power such that $q\equiv 1\,(\mod{4})$.  
\end{itemize}
In particular, there exists a Hadamard matrix of order $4(2n+1)$ if $n$ satisfies any of the above conditions. 
\end{theorem}
In this paper, we obtain a new series of difference families of type $H_4^*$ with parameters 
$(2q^2;q^2,q^2,$ $q^2,q^2-1;2q^2-2)$ where $q$ is a prime power congruent to 3 modulo 8 satisfying some extra condition. The construction uses $8^{\rm th}$ cyclotomic  classes of $\F_{q^2}$ and ``half lines" in $\AG(2,q)$.  In \cite{LMS06,XL91,XL96}, the main idea is to construct difference families of type $H$ in the group $(\F_{q^2},+)$. Our approach here is analogous to that of \cite{LM18}; the main difference here is the usage of Paley type partial difference sets. 
The following are our main results.  
\begin{theorem}\label{thm:mainth}
Let $q$ be a prime power of the form $q=12c^2+4c+3$ with 
$c$ an arbitrary integer, and let $n=q^2$.  Then there exists a difference family with parameters $(2n;n,n,n,n-1;2n-2)$ in $\Z_{2}\times (\F_{q^2},+)$. 
\end{theorem}\
By plugging the group invariant $(1,-1)$ matrices obtained from the blocks of the difference family in Theorem~\ref{thm:mainth} into the Wallis-Whiteman array, we immediately obtain the following: 
\begin{theorem}\label{thm:mainthm}
Let $q$ be a prime power of the form $q=12c^2+4c+3$ with 
$c$ an arbitrary integer, and let $n=q^2$.  Then there exists a  
Hadamard matrix of order $4(2n+1)$. 
\end{theorem}
We remark that there are $386$ prime powers of the form $q=12c^2+4c+3<10^7$ while there are $166181$ prime powers $q<10^7$ such that $q\equiv 3\,(\mod{8})$. The first $58$ prime powers of  the form $q=12c^2+4c+3<10^5$ are listed below: 
\begin{align}
&3, 11, 19, 43, 59, 179, 211, 283, 563, 619, 739, 1163, 1499, 1979, 2083, 2411, 3011,\nonumber\\
&3539, 4259, 4723, 7603, 8011, 8219, 10211, 11411, 12163, 14011, 14563, 14843,\nonumber\\ 
&17483, 20011, 23059, 25579, 26699, 28619, 29803, 30203, 33923, 36083, 36523, \label{eq:list}\\
& 41539, 49411, 54139, 55219, 55763, 59083, 60779, 63659, 65419, 69011, 70843, \nonumber\\
& 75211, 80363, 81019, 82339, 83003, 88411, 93283.\nonumber
\end{align}

%%%%%%%%%%%%%%%%%%%%%%%%%%%%%%%%%%%%%%%%%%%%%%%%%%%%%%%%%%%%%
\section{The construction}

We first fix our notation.  Let $q$ be a prime power such that $q\equiv 3\,(\mod{4})$.  Let $\omega$ be a primitive element of $\F_{q^2}$ and 
let $0_{\F_{q^2}}$ denote the zero of $\F_{q^2}$. For any fixed positive integer $N$ dividing $q^2-1$, define $C_i^{(N,q^2)}=\omega^i \langle \omega^N\rangle$, $i=0,1,\ldots,N-1$, 
called the {\it $N^{\rm th}$ cyclotomic classes} of $\F_{q^2}$. Furthermore, define   
\begin{align*}
H_i=&\,C_i^{(2(q+1),q^2)}, \, \, i=0,1,\ldots,2q+1,\\
L_i=&\,C_i^{(q+1,q^2)}, \, \, i=0,1,\ldots,q,\\
S_i=&\,C_i^{(q+1,q^2)}\cup \{0\}, \, \, i=0,1,\ldots,q, \\
 D_i= & C_i^{(4,q^2)}\cup C_{i+1}^{(4,q^2)},  \, \, i=0,1,\ldots,3.
\end{align*}
Note that each $S_i$ is a line through the origin of $\AG(2,q)$; for this reason the $H_i$'s are called half lines \cite{xiang}. In the group ring $\Z[(\F_{q^2},+)]$, we have
\begin{equation}\label{eq:SS}
\mbox{$S_iS_j=\F_{q^2}$ for $i\not=j$ and $S_i^2=qS_i$ for all $i$.} 
\end{equation}

\begin{lemma}
\label{lem:pds}
%For any prime power $q\equiv 1\,(\mod{4})$, the set $D=C_0^{(2,q)}$ of nonzero squares satisfies that $DD^{(-1)}=\frac{q-5}{4}D+\frac{q-1}{4}\overline{D}+\frac{q-1}{2[0_{\F_q}]$. The set  $D$ is called the Paley partial difference set. On the other hand, 
For $i=0,1,2,3$, $D_i$ is a Paley type partial difference set in $(\F_{q^2},+)$. In particular, 
\[ D_iD_i^{(-1)}=\frac{q^2-5}{4}D_i+\frac{q^2-1}{4}\overline{D}_i+\frac{q^2-1}{2}\cdot 0_{\F_{q^2}}. \]
\end{lemma}

For a proof of  Lemma~\ref{lem:pds}, we refer the reader to  \cite[p.~216]{peisert}. The strongly regular Cayley graph, ${\rm Cay}(\F_{q^2}, D_0)$, is often called a Peisert graph.

Our objective is to construct difference families with parameters $(2q^2;q^2,q^2-1,q^2,q^2;2q^2-2)$ 
in $\Z_2\times \F_{q^2}$. So we need to find four blocks $B_0, B_1, B_2, B_3$ with $|B_i|=q^2$, $i=0,2,3$, and $|B_1|=q^2-1$, in $\Z_2\times \F_{q^2}$ such that 
\[ \sum_{i=0}^{3} B_i B_i^{(-1)}= (2q^2-2) (\Z_2\times \F_{q^2}) + (2q^2+1) \cdot (0,0_{\F_{q^2}}).\]

To construct the first two blocks, we make use of the Paley type partial difference sets $D_0$ and $D_2$ defined above.  Note that $D_2=\omega^2 D_0$ and $D_2=\overline{D_0}$. In $\Z_{2}\times \F_{q^2}$, we set 
\begin{align*}
B_0&\,=(\{0\}\times D_0)\cup (\{1\}\times (\F_{q^2}\setminus D_0)), \\
B_1&\,=(\{0\}\times D_2)\cup (\{1\}\times D_2).
\end{align*}

Then $|B_0|=q^2$ and $|B_1|=q^2-1$. 

\begin{proposition}\label{prop:b0b1}
With $B_0, B_1$ defined as above, we have
\begin{equation}\label{eq:firstB}
\sum_{i=0,1}B_iB_i^{(-1)}= 
\{0\}\times \Big(
(q^2-2)\F_{q^2}^\ast+(2q^2-1)\cdot 0_{\F_{q^2}}\Big)+\{1\}\times \Big(2D_0-2D_2
+(q^2-1)\F_{q^2}\Big).
\end{equation}
\end{proposition}
\proof
It is clear that 
\begin{align}
\sum_{i=0,1}B_iB_i^{(-1)}=&\, 
\{0\}\times \Big(2\sum_{i=0,2}D_iD_i^{(-1)}+(q^2-2|D_0|)\F_{q^2}\Big)\nonumber\\
&\, \, +\{1\}\times \Big(-2D_0D_0^{(-1)}+2D_2D_2^{(-1)}
+2|D_0|\F_{q^2}\Big). \label{eq:to1}
\end{align}
By Lemma~\ref{lem:pds}, we have 
\begin{align}
\sum_{i=0,2}D_iD_i^{(-1)}=&\,\frac{q^2-5}{4}(D_0+D_2)+\frac{q^2-1}{4}
(\overline{D}_0+\overline{D}_2)+(q^2-1)\cdot 0_{\F_{q^2}}\nonumber\\
=&\, \frac{q^2-3}{2}\F_{q^2}^\ast+(q^2-1)\cdot 0_{\F_{q^2}},   \label{eq:D1}
\end{align}
and 
\begin{align}
-D_0D_0^{(-1)}+D_2D_2^{(-1)}=&\,-(\frac{q^2-5}{4}D_0+\frac{q^2-1}{4}D_2)+
(\frac{q^2-5}{4}D_2+\frac{q^2-1}{4}D_0)\nonumber\\
=&\,D_0-D_2. \label{eq:D2}
\end{align}
It is now straight forward to obtain \eqref{eq:firstB} from \eqref{eq:to1}, \eqref{eq:D1} and \eqref{eq:D2}. 
%\[
%\sum_{i=0,1}B_iB_i^{(-1)}= 
%\{0\}\times \Big(
%(q^2-2)\F_{q^2}^\ast+(2q^2-1)[0_{\F_{q^2}}]\Big)+\{1\}\times \Big(2D_0-2D_2
%+2(q^2-1)\F_{q^2}\Big).
%\]
%This completes the proof of the proposition. 
\qed
\vspace{0.3cm}

To construct the remaining blocks of the desired difference family, we need difference families of type $H$ in $\F_{q^2}$ that satisfy certain conditions. 

\begin{proposition}\label{prop:b0b2}
Suppose  ${\mathcal E}=\{E_i:i=0,1,2,3\}$ is a difference family of type $H$ in $\F_{q^2}$ such that $|E_0|=|E_1|=|E_2|=|E_3|=(q^2-q)/2$ and 
\begin{equation}\label{eq:cand_ho}
E_0E_1^{(-1)}+E_1E_0^{(-1)}+E_2E_3^{(-1)}+E_3E_2^{(-1)}=(q-1)^2\F_{q^2}+2D_0
-2D_2.  
\end{equation}
Let $B_0, B_1$ be defined as above and set
\begin{align*}
B_2&\,=(\{0\}\times E_0)\cup (\{1\}\times (\F_{q^2}\setminus E_1)), \\
B_3&\,=(\{0\}\times E_2)\cup (\{1\}\times (\F_{q^2}\setminus E_3)).
\end{align*}
Then $\{B_0, B_1, B_2,B_3\}$ is a difference family with parameters $(2q^2;q^2,q^2-1,q^2,q^2;2q^2-2)$ 
in $\Z_2\times \F_{q^2}$. 
\end{proposition}
\proof  First of all, we have $|B_2|=q^2+|E_0|-|E_1|=q^2$ and $|B_3|=q^2+|E_2|-|E_3|=q^2$. In view of \eqref{eq:firstB}, it suffices to show that 
\[
\sum_{i=2,3}B_iB_i^{(-1)}= 
\{0\}\times (
2q^2\cdot 0_{\F_{q^2}}+q^2 \F_{q^2}^\ast)+\{1\}\times ((q^2-1)\F_{q^2}-2D_0+2D_2). 
\]
It is clear that 
\begin{align}
\sum_{i=2,3}B_iB_i^{(-1)}=&\, 
\{0\}\times \Big(\sum_{i=0}^3E_iE_i^{(-1)}+2(q^2-|E_1|-|E_3|)\F_{q^2}\Big)\nonumber\\
&\, \, +\{1\}\times (-E_0E_1^{(-1)}-E_1E_0^{(-1)}
-E_2E_3^{(-1)}-E_3E_2^{(-1)}
+2(|E_0|+|E_2|)\F_{q^2}). \label{eq:to2}
\end{align}
Since $\{E_i:i=0,1,2,3\}$ is a difference family of type $H$ and 
$|E_1|+|E_3|=q^2-q$, we have 
\begin{equation}\label{eq:e1}
\sum_{i=0}^3E_iE_i^{(-1)}+2(q^2-|E_1|-|E_3|)\F_{q^2}\\
=2q^2 \cdot 0_{\F_{q^2}}+q^2 \F_{q^2}^\ast.  
\end{equation}
On the other hand, 
by the assumption~\eqref{eq:cand_ho} and $|E_0|+|E_2|=q^2-q$, we have 
\begin{equation}
-E_0E_1^{(-1)}-E_1E_0^{(-1)}
-E_2E_3^{(-1)}-E_3E_2^{(-1)}
+2(|E_0|+|E_2|)\F_{q^2}
=(q^2-1)\F_{q^2}-2D_0+2D_2.\label{eq:e2}
\end{equation}
The proposition now follows from \eqref{eq:to2}, \eqref{eq:e1}, and \eqref{eq:e2}.  
\qed

\medskip

%From now on, we use the following notations. 

To construct difference families of type $H$ in $\F_{q^2}$ satisfying the conditions in Proposition~\ref{prop:b0b2}, it is then natural to consider those constructed in \cite{LMS06}. 
\begin{lemma}{\em (\cite[Lemma~4 and Corollary~5]{LMS06})}
\label{lem:scl}
Let $q\equiv 3\,(\mod{4})$ be a prime power and let $e$ be the exact power of $2$ dividing $q+1$. 
Let $\alpha<e$ be an odd number and set $\beta=\frac{qe-\alpha(q+1)}{2e}$. Let ${\bf A}\subseteq \{0,1,\ldots,2e-1\}$ and ${\bf B}_0,\ldots,{\bf B}_{e-1}\subseteq \{0,1,\ldots,q\}$ with $|{\bf A}|=\alpha$, 
$|{\bf B}_0|=\cdots=|{\bf B}_{e-1}|=\beta$ such that $b\not\equiv a\,(\mod{e})$ for all $a\in {\bf A}$ and $b\in \bigcup_{r=0}^{e-1}{\bf B}_r$. Set \begin{align*}
H&=\bigcup_{i\in {\bf A}}C_{i}^{(2e,q^2)}\\
M_i&=\bigcup_{j\in {\bf B}_i}L_j, \, \, i=0,1,\ldots,e-1\\
{\bf D}_i&=\omega^i(H\cup M_i), \, \, i=0,1,\ldots,e-1. 
\end{align*}
Then $|{\bf D}_i|=\frac{q(q-1)}{2}$ for  $i=0,1,\ldots,e-1$, and $\{{\bf D}_i: i=0,1,\ldots,e-1\}$ forms a
difference family in $(\F_{q^2},+)$ with $\lambda=\frac{eq(q-2)}{4}$.  
\end{lemma}

We now assume that $q$ is a prime power and $q=8m+3$ for some positive integer $m$. In view of Lemma~\ref{lem:scl}, we need a set ${\bf A}\subseteq \{0,1,\ldots,7\}$ with $|{\bf A}|=3$, and four subsets ${\bf B}_i$, $i=0,1,2,3$, of $\{0,\ldots, q\}$, each of size $m$, satisfying certain conditions. 

First, we require $I \cap \{x+4\,(\mod{8}):x \in I\}=\emptyset$. Since $|I|=3$, the condition $I \cap \{x+4\,(\mod{8}):x \in I\}=\emptyset$ simply means that $I$ contains exactly one odd or exactly one even element, say, $y\in I$. (Note that such an $I$ clearly exists, for example, take $I=\{0,1,3\}$; and in this case $y=0$.)  Next, we define two $m$-subsets of $\{0,1,\ldots ,q\}$: 
\[ J_1=\{y+2+4i\,(\mod{q+1}):i\in \{0,1,\ldots,m-1\}\} \mbox{ and } \]
\[ J_2=\{y+4i\,(\mod{q+1}):i\in \{0,1,\ldots,m-1\}\}.\] %Then, these sets  $I_1,I_2,J_1,J_2$ satisfy the conditions (1)--(4) of Lemma~\ref{lem:basic1}. 

Now, using the notation in Lemma~\ref{lem:scl}, we set $e=4$, $\alpha=3$ and $\beta=m$. 
Let ${\bf A}=I$, ${\bf B}_0={\bf B}_1=J_1$,
\[{\bf B}_2={\bf B}_3=\{y-2+4i\,(\mod{q+1}):i\in \{0,1,\ldots,m-1\}\}.\]
It is then straight forward to check that the conditions in Lemma~\ref{lem:scl} are all satisfied. Therefore we obtain a difference family $\{{\bf D}_i:i=0,1,2,3\}$. However, for our purpose, we need to set $E_0={\bf D}_0, E_1={\bf D}_2, E_2={\bf D}_1$ and $E_3={\bf D}_3$. In terms of $I, J_1, J_2$, we have the following:  
\begin{align}
E_0=&\,\Big(\bigcup_{i\in I}C_i^{(8,q^2)}\Big)\cup \Big(\bigcup_{i\in J_1}L_i\Big), \, E_1=\Big(\bigcup_{i\in I}C_{i+2}^{(8,q^2)}\Big)\cup \Big(\bigcup_{i\in J_2}L_i\Big), \label{eq:blocks}\\
E_2=&\,\Big(\bigcup_{i\in I}C_{i+1}^{(8,q^2)}\Big)\cup \Big(\bigcup_{i\in J_1}L_{i+1}\Big), \, 
E_3=\Big(\bigcup_{i\in I}C_{i+3}^{(8,q^2)}\Big)\cup \Big(\bigcup_{i\in J_2}L_{i+1}\Big). \nonumber
\end{align}
By Lemma~\ref{lem:scl}, $\{ E_i, i=0,1,2,3\}$ is a difference family of type $H$  in $(\F_{q^2},+)$. Furthermore, 
$|E_i|=\frac{q^2-q}{2}$ for $i=0,1,2,3$. It therefore remains to show the following:
%Hence, we have 
%\begin{equation}\label{eq:Es}
%\sum_{i=0}^3E_iE_i^{(-1)}=2(q^2-q)[0_{\F_{q^2}}]+q(q-2)\F_{q^2}^\ast. 
%\end{equation}
%Hereafter, we co

\begin{theorem}\label{prop:b0b13}
The $E_i$'s defined in (2.10) satisfy the equation (2.6). In particular, there is a difference family with parameters $(2q^2;q^2,q^2,q^2,q^2-1;2q^2-2)$ 
in $\Z_2\times (\F_{q^2},+)$. 
\end{theorem}

\section{Proof of Theorem ~\ref{prop:b0b13}}

To prove Theorem ~\ref{prop:b0b13}, we need to compute 
$E_0E_1^{(-1)}+E_0E_1^{(-1)}+E_2E_3^{(-1)}+E_3E_2^{(-1)}$. As in the case of Lemma 4 in \cite{LMS06}, it will make the computations easier if we write each $E_i$ in a different form (i.e., as a union of $H_i$'s and $L_j$'s).  Recall that $q=8m+3$ is a prime power. We define 
\[  I_1=\{x+8i\,(\mod{2(q+1)}) :x\in I, i\in \{0,1,\ldots,2m\}\} \mbox{ and } I_2=I_1+2.\]
Here we use the notation $K+1=\{x+1: x\in K\}$. Note that $|I_1|=|I_2|=3(q+1)/4$. Recall that 
\[ J_1=\{y+2+4i\,(\mod{q+1}):i\in \{0,1,\ldots,m-1\}\} \mbox{ and }  J_2=J_1-2.\] 
We write 
\[ E_0=\sum_{i\in I_1}H_i+\sum_{i\in J_1}L_i \mbox{ and } 
E_1=\sum_{i\in I_2}H_i+\sum_{i\in J_2}L_i, \]
\[ E_2=\sum_{i\in I_1+1}H_i+\sum_{i\in J_1+1}L_i \mbox{ and } 
E_3=\sum_{i\in I_2+1}H_i+\sum_{i\in J_2+1}L_i. \]
 Observe that the following conditions are satisifed: 
\begin{itemize}
\item[(1)] Since $I \cap \{x+4\,(\mod{8}):x \in I\}=\emptyset$, we have $I_i\cap \{h+(q+1)\,(\mod{2(q+1)}):h\in I_i\}=\emptyset$ for $i=1,2$, 
\item[(2)] $a\not\equiv b\,(\mod{q+1})$ for all $a\in I_i,b\in J_i$, $i=1,2$,
\item[(3)] $|I_1|+2|J_1|=|I_2|+2|J_2|=q$, 
\item[(4)] $J_1\subseteq I_2'\cup J_2$ and  $J_2\subseteq I_1'\cup J_1$,  %and $J_1\cap J_2=\emptyset$, 
where $I_j'=\{i\,(\mod{q+1}):i \in I_j\}$ for $j=1,2$. 
\end{itemize}

\begin{lemma} \label{lem:basic1}
In the group ring $\Z[(\F_{q^2},+)]$, $E_0E_1^{(-1)}+E_1E_0^{(-1)}=$ 
\begin{equation}\label{eq:x1x2}
\sum_{i\in I_1}\sum_{j\in I_2}H_iH_j^{(-1)}+
\sum_{i\in I_2}\sum_{j\in I_1}H_iH_j^{(-1)}
+\lambda_1 \cdot 0_{\F_{q^2}}
+\lambda_2\F_{q^2} -|J_2|\sum_{i\in I_1'}S_i-|J_1|\sum_{i\in I_2'}S_i, 
\end{equation}
where $\lambda_1=|I_1||J_2|+|I_2||J_1|+2|J_1||J_2|$ and 
$\lambda_2=|I_1'||J_2|+|I_2'||J_1|+2|J_1||J_2|-|I_1'\cap J_2|-|I_2'\cap J_1|-2|J_1\cap J_2|$. 
\end{lemma}
\proof Note that $L_i^{(-1)}=L_i$. We first expand the expression  $E_0E_1^{(-1)}+E_0E_1^{(-1)}$ and obtain the following:
\[E_0E_1^{(-1)}+E_1E_0^{(-1)}=\sum_{i\in I_1}\sum_{j\in I_2}H_iH_j^{(-1)}+
\sum_{i\in I_2}\sum_{j\in I_1}H_iH_j^{(-1)}+Y \mbox{ where } \]
\[ Y= \sum_{i\in I_1}(H_i+H_i^{(-1)})\sum_{i\in J_2}L_i +\sum_{i\in I_2}(H_i+H_i^{(-1)})\sum_{i\in J_1}L_i + 2\sum_{i\in J_1}L_i\sum_{i\in J_2}L_i. 
\]

Note that $S_i=L_i+0_{\F_{q^2}}$. So, we may replace each $L_i$ by $S_i-0_{\F_{q^2}}$ in the above sum and we get
\begin{align}
Y=&\, -|J_2|\sum_{i\in I_1}(H_i+H_i^{(-1)}) -|J_1|\sum_{i\in I_2}(H_i+H_i^{(-1)}) \nonumber\\
 & + \sum_{i\in I_1}\sum_{j\in J_2}S_j(H_i+H_i^{(-1)})+
\sum_{i\in I_2}\sum_{j\in J_1}S_j(H_i+H_i^{(-1)})\nonumber\\
&+2\sum_{i\in J_1}\sum_{j\in J_2}S_iS_j
-2|J_1|\sum_{j\in J_2}S_i-2|J_2|\sum_{j\in J_1}S_i+2|J_1||J_2|\cdot 0_{\F_{q^2}}.
\nonumber
\end{align}

Observe that  $H_i+H_{q+1+i}=S_i- 0_{\F_{q^2}}$ for $i=0,1,\ldots,q$ and (3) holds. 
Also note that 
\[ \sum_{i\in I_1}\sum_{j\in J_2}S_j(H_i+H_i^{(-1)})= \sum_{i\in I_1'}\sum_{j\in J_2}S_jS_i-|J_2|\sum_{i\in I_1'}S_i \mbox{ and } \]
\[ \sum_{i\in I_2}\sum_{j\in J_1}S_j(H_i+H_i^{(-1)})= \sum_{i\in I_2'}\sum_{j\in J_1}S_jS_i-|J_1|\sum_{i\in I_2'}S_i.\]
We then have 
\begin{align}
Y =&\, \lambda_1 \cdot 0_{\F_{q^2}} -|J_2|\sum_{i\in I_1'}S_i -|J_1|\sum_{i\in I_2'}S_i - q\sum_{i\in J_2}S_i -q\sum_{i\in J_1}S_i \nonumber\\
 & 
+\sum_{i\in I_1'}\sum_{j\in J_2}S_jS_i+
\sum_{i\in I_2'}\sum_{j\in J_1}S_iS_j
+2\sum_{i\in J_1}\sum_{j\in J_2}S_iS_j.\label{eq:compx1x2}
\end{align}

On the other hand, $S_iS_j=\F_{q^2}$ whenever $i\neq j$. Therefore, 
by the conditions~(2) and (4), for distinct $u,v$ in $\{1,2\}$,
\begin{equation}\label{eq:SS2}
\sum_{i\in I_u'}\sum_{j\in J_v}S_jS_i=q\sum_{i\in (I_u'\cap J_v)}S_i+(|I_u'|\cdot |J_v|-|I_u\cap J_v|)\F_{q^2}  \mbox{ and }
\end{equation} 
\begin{equation}\label{eq:SS3}
\sum_{i\in J_1}\sum_{j\in J_2}S_iS_j=q\sum_{i\in (J_1\cap J_2)}S_i+(|J_1|\cdot |J_2|-|J_1\cap J_2|)\F_{q^2}.
\end{equation}
\eqref{eq:x1x2} now follows easily from \eqref{eq:compx1x2},  \eqref{eq:SS2} and \eqref{eq:SS3}.
\qed

\medskip

Now, replace $I_i$ with $I_i+1$, and $J_i$ with $J_i+1$ in the argument above and observe that condition (2), (3) and (4) still hold. We immediately get the following:

\begin{lemma} \label{lem:basic2}
In the group ring $\Z[(\F_{q^2},+)]$, $E_2E_3^{(-1)}+E_3E_2^{(-1)}=$ 
\begin{equation}\label{eq:x1x3}
\sum_{i\in I_1+1}\sum_{j\in I_2+1}H_iH_j^{(-1)}+
\sum_{i\in I_2+1}\sum_{j\in I_1+1}H_iH_j^{(-1)}
+\lambda_1 \cdot 0_{\F_{q^2}}
+\lambda_2\F_{q^2} -|J_2|\sum_{i\in I_1'}S_i-|J_1|\sum_{i\in I_2'}S_i, 
\end{equation}
where $\lambda_1=|I_1||J_2|+|I_2||J_1|+2|J_1||J_2|$ and 
$\lambda_2=|I_1'||J_2|+|I_2'||J_1|+2|J_1||J_2|-|I_1'\cap J_2|-|I_2'\cap J_1|-2|J_1\cap J_2|$. 
\end{lemma}

\medskip

\begin{lemma}\label{lem:eout}
Let $E_i$, $i=0,1,2,3$, be defined as in \eqref{eq:blocks}. Recall that $q=8m+3$. Then, we have
\[ E_0E_1^{(-1)}+E_0E_1^{(-1)}+E_2E_3^{(-1)}+E_3E_2^{(-1)}=\sum_{h=0,1}\sum_{i,j\in I}\Big(C_{i+h}^{(8,q^2)}{C_{j+2+h}^{(8,q^2)}}^{(-1)}+
C_{i+2+h}^{(8,q^2)}{C_{j+h}^{(8,q^2)}}^{(-1)}\Big)+Z \]
where $Z=8m(4m+1)\cdot 0_{\F_{q^2}}+m(28m+5)\F_{q^2}^\ast$. 
\end{lemma}
\proof  Applying Lemmas \ref{lem:basic1} and \ref{lem:basic2}, we obtain 
\begin{align}
&\,E_0E_1^{(-1)}+E_0E_1^{(-1)}+E_2E_3^{(-1)}+E_3E_2^{(-1)}\nonumber\\
=&\, \sum_{h=0,1}\sum_{i,j\in I}\Big(C_{i+h}^{(8,q^2)}{C_{j+2+h}^{(8,q^2)}}^{(-1)}+
C_{i+2+h}^{(8,q^2)}{C_{j+h}^{(8,q^2)}}^{(-1)}\Big)
+2\lambda_1 0_{\F_{q^2}}
+2\lambda_2\F_{q^2}\label{eq:CC}
\\
&\hspace{1.3cm} -|J_2|\sum_{i\in I_1'}(S_i+S_{i+1})-|J_1|\sum_{i\in I_2'}(S_i+S_{i+1}). \nonumber
\end{align}
Since $|I_1'|=|I_2'|=6m+3$,  we have 
\[
\sum_{i\in I_1'}(S_i+S_{i+1})=2(6m+3)\cdot 0_{\F_{q^2}}+\sum_{i\in I}
\Big(C_i^{(4,q^2)}+C_{i+1}^{(4,q^2)}\Big)
\]
and 
\[
\sum_{i\in I_2'}(S_i+S_{i+1})=2(6m+3)\cdot 0_{\F_{q^2}}+\sum_{i\in I}
\Big(C_{i+2}^{(4,q^2)}+C_{i+3}^{(4,q^2)}\Big). 
\]
Note that $\sum_{j=0}^3 C_{i+j}^{(4,q^2)}=\F_{q^2}^\ast$, $|I|=3$ and $|J_1|=|J_2|=m$. Hence, 
\begin{equation}\label{eq:CC1}
-|J_2|\sum_{i\in I_1'}(S_i+S_{i+1})-|J_1|\sum_{i\in I_2'}(S_i+S_{i+1})=
-12m(2m+1)\cdot 0_{\F_{q^2}}-3m\F_{q^2}^\ast. 
\end{equation}
Furthermore, it is clear that 
\begin{equation}\label{eq:CC2}
\lambda_1=2m(7m+3) \mbox{\, \, and \, } 
\lambda_2=2m(7m+2). 
\end{equation}
Our lemma now follows from \eqref{eq:CC} with \eqref{eq:CC1} and \eqref{eq:CC2}. 
\qed
\vspace{0.3cm}

%\begin{proposition}\label{prop:b2b3}
%Define 
%\begin{align*}
%B_2=&\,\{0\}\times E_0 \cup \{1\}\times (\F_{q^2}\setminus E_1), \\
%B_3=&\,\{0\}\times E_2\cup \{1\}\times (\F_{q^2}\setminus E_3).  
%\end{align*}
%Assume that 
%\begin{equation}\label{eq:cccc}
%\sum_{h=0,1}\sum_{i,j\in I}\Big(C_{i+h}^{(8,q^2)}{C_{j+2+h}^{(8,q^2)}}^{(-1)}+
%C_{i+2+h}^{(8,q^2)}{C_{j+h}^{(8,q^2)}}^{(-1)}\Big)=
%\frac{q^2-1}{2}[0_{\F_{q^2}}]+\frac{9q^2-17}{16}\F_{q^2}^\ast+2D_0-2D_2. 
%\end{equation}
%Then, it holds that 
%\begin{equation}\label{eq:cand2}
%\sum_{i=2,3}B_iB_i^{(-1)}= 
%\{0\}\times (
%2q^2[0_{\F_{q^2}}]+q^2 \F_{q^2}^\ast)+\{1\}\times ((q^2-1)\F_{q^2}-2D_0+2D_2). 
%\end{equation}
%\end{proposition}
%\proof
%It is clear that 
%\begin{align}
%\sum_{i=2,3}B_iB_i^{(-1)}=&\, 
%\{0\}\times \Big(\sum_{i=0}^3E_iE_i^{(-1)}+2(q^2-|E_1|-|E_3|)\F_{q^2}\Big)\nonumber\\
%&\, \, +\{1\}\times (-E_0E_1^{(-1)}-E_1E_0^{(-1)}
%-E_2E_3^{(-1)}-E_3E_2^{(-1)}
%+2(|E_0|+|E_2|)\F_{q^2}). \label{eq:to2}
%\end{align}
%By \eqref{eq:Es}, we have 
%\begin{equation}\label{eq:e1}
%\sum_{i=0}^3E_iE_i^{(-1)}+2(q^2-|E_1|-|E_3|)\F_{q^2}\\
%=2q^2[0_{\F_{q^2}}]+q^2 \F_{q^2}^\ast.  
%\end{equation}
%On the other hand, 
%by Lemma~\ref{lem:eout} with our assumption, we have 
%\begin{equation}
%-E_0E_1^{(-1)}-E_1E_0^{(-1)}
%-E_2E_3^{(-1)}-E_3E_2^{(-1)}
%+2(|E_0|+|E_2|)\F_{q^2}
%=(q^2-1)\F_{q^2}-2D_0+2D_2.\label{eq:e2}
%\end{equation}
%Hence, continuing from \eqref{eq:to2} with \eqref{eq:e1} and \eqref{eq:e2}, %we obtain 
%the equation~\eqref{eq:cand2}. 
%\qed

To finish our proof, we need to evaluate  $C_i^{(8,q^2)}{C_j^{(8,q^2)}}^{(-1)}$. The coefficient $c_x$ of $x\in \F_{q^2}$ in $C_i^{(8,q^2)}{C_j^{(8,q^2)}}^{(-1)}$ is $\big|(C_j^{(8,q^2)}+x)\cap C_i^{(8,q^2)}\big|$. If $x\in C_h^{(8,q^2)}$,  
it is clear that $c_x=\big|(C_{j-h}^{(8,q^2)}+1)\cap C_{i-h}^{(8,q^2)}\big|$. 
The numbers $(i,j)_{N}=\big|(C_{i}^{(N,q^2)}+1)\cap C_{j}^{(N,q^2)}\big|$, $i,j=0,1,\ldots,N-1$, are called {\it $N^{\rm th}$ cyclotomic numbers}. In our case, $q\equiv 3\,(\mod{8})$ is a prime power.
In view of \cite[Lemma~30]{S67}, we obtain the following:

\begin{proposition} \label{prop:cyclo8}
Let $q\equiv 3\,(\mod{8})$ be a prime power. 
Then the cyclotomic numbers $(i,j)_{8}$, $i,j=0,1,\ldots,7$, in $\F_{q^2}$ are determined by Table~\ref{Tab1} and the relations: 
\begin{align*}
&64n_1=q^2-15+2q, \, 64n_2=q^2+1-2q-4a,\, 
64n_3=q^2+1-6q+8a, \\
&64n_4=q^2+1+18q, \, 64n_5=q^2-7-2q+4a, \, 
64n_6=q^2+1+6q+4a+16b, \\ 
&64n_7=q^2+1+6q+4a-16b, \, 
64n_8=q^2-7+2q-8a, 
\end{align*}
where $a,b$ are specified by the unique proper representation of 
$q^2=a^2+2b^2$ with $a\equiv 1\,(\mod{4})$. Note that there is no restriction on the 
sign of $b$.

{\small 
\begin{table}[h]
\caption{Cyclotomic numbers of order $8$: 
the $(i,j)$-entry is  
$(i,j)_8$.}
\label{Tab1}
\vspace{-0.4cm}
$$
\begin{array}{|c||c|c|c|c|c|c|c|c|}
\hline
&0&1&2&3&4&5&6&7\\
\hline
\hline 
0&n_1&n_2&n_3&n_2&n_4&n_2&n_3&n_2\\
\hline
1&n_5&n_5&n_6&n_2&n_2&n_2&n_2&n_7\\
\hline
2&n_8&n_2&n_8&n_7&n_3&n_2&n_3&n_6\\
\hline
3&n_5&n_2&n_2&n_5&n_2&n_7&n_6&n_2\\
\hline
4&n_1&n_5&n_8&n_5&n_1&n_5&n_8&n_5\\
\hline
5&n_5&n_2&n_7&n_6&n_2&n_5&n_2&n_2\\
\hline
6&n_8&n_7&n_3&n_2&n_3&n_6&n_8&n_2\\
\hline
7&n_5&n_6&n_2&n_2&n_2&n_2&n_7&n_5\\
\hline
\end{array}
$$
\end{table}}
\end{proposition}

\begin{theorem}\label{thm:b2b3}
Suppose $q^2=a^2+2b^2$ is the unique proper representation with $a\equiv 1\,(\mod{4})$.  Theorem~\ref{prop:b0b13} holds if either of the following conditions is satisfied.
\begin{itemize}
\item [(a)] $I=\{0,2,3\}$ and $3q=a+ 4b+16$.
\item [(b)] $I=\{0,2,7\}$ and  $3q=a- 4b+16$. 
\end{itemize}
\end{theorem}
\proof By Lemma~\ref{lem:eout}, it is sufficient to show the following: 
\begin{equation}\label{eq:cccc}
U:=\sum_{h=0,1}\sum_{i,j\in I}\Big(C_{i+h}^{(8,q^2)}{C_{j+2+h}^{(8,q^2)}}^{(-1)}+
C_{i+2+h}^{(8,q^2)}{C_{j+h}^{(8,q^2)}}^{(-1)}\Big)=
\frac{q^2-1}{2}\cdot 0_{\F_{q^2}}+\frac{9q^2-17}{16}\F_{q^2}^\ast+2D_0-2D_2. 
\end{equation}

We give a proof only in the case where $3q=a+4b+16$. The proof for 
the case where $3q=a-4b+16$ is similar. 

Define $D=\bigcup_{i\in I}C_i^{(8,q^2)}$, and 
let $c_x$ denote the 
coefficient  of $x\in \F_{q^2}$ in $U$. 
%\[
%\sum_{h=0,1}\sum_{i,j\in I}\Big(C_{i+h}^{(8,q^2)}{C_{j+2+h}^{(8,q^2)}}^{(-1)}+
%C_{i+2+h}^{(8,q^2)}{C_{j+h}^{(8,q^2)}}^{(-1)}\Big). 
%\] 
To show that  $c_0=\frac{q^2-1}{2}$, it is sufficient to check number of pairs $(i,j)\in I\times I$ such that
$i\equiv j+2\,(\mod{8})$ and $i+2\equiv j\,(\mod{8})$. Clearly, the solution is $(2,0)$ and $(0,2)$ in each case respectively. Therefore, $c_0=2\times 2 \times \frac{q^2-1}{8}=\frac{q^2-1}{2}$. 
 To prove that \eqref{eq:cccc} holds, it is enough to 
see that $c_{1}=c_{\omega}=c_{\omega^2}+4=c_{\omega^3}+4$ since 
$c_{w^i}=c_{w^{i+4j}}$ for all $i,j$.  
On the other hand, 
$c_x$ for $x\in \F_{q^2}^\ast$ 
is given by 
\[
c_x=|D\cap (\omega^2 D+x)|+|\omega^2 D\cap (D+x)|+
|\omega D\cap (\omega^3 D+x)|+|\omega^3 D\cap (\omega D+x)|. 
\]
Hence, the system of equations $c_{1}=c_{\omega}=c_{\omega^2}+4=c_{\omega^3}+4$ is reformulated as 
\begin{align*}
&|D\cap (\omega^2 D+1)|+|\omega^2 D\cap  (D+1)|
+|\omega D\cap (\omega^3 D+1)|+|\omega^3 D\cap  (\omega D+1)|\\
=&\,|\omega^{-1} D\cap (\omega D+1)|+|\omega D\cap  (\omega^{-1} D+1)|+| D\cap (\omega^2 D+1)|+|\omega^2 D\cap  (D+1)|\\
=&\,|\omega^{-2} D\cap (D+1)|+|D\cap  (\omega^{-2} D+1)|+| \omega^{-1}D\cap (\omega D+1)|+|\omega D\cap  (\omega^{-1}D+1)|+4\\
=&\,|\omega^{-3} D\cap (\omega^{-1}D+1)|+|\omega^{-1}D\cap  (\omega^{-3} D+1)|+| \omega^{-2}D\cap (D+1)|+| D\cap  (\omega^{-2}D+1)|+4. 
\end{align*}
Noting that $|\omega D\cap (\omega^3 D+1)|+|\omega^3 D\cap (\omega D+1)|=|\omega^{-3} D\cap (\omega^{-1} D+1)|+|\omega^{-1} D\cap (\omega^{-3} D+1)|$, 
the equations above are reduced to 
\begin{equation}\label{eq:redu1}
|\omega D\cap (\omega^3 D+1)|+|\omega^3 D\cap (\omega D+1)|=
|\omega^{-1} D\cap (\omega D+1)|+|\omega D\cap  (\omega^{-1} D+1)|
\end{equation}
and 
\begin{equation}\label{eq:redu2}
|D\cap (\omega^2 D+1)|+|\omega^2 D\cap  (D+1)|=
|\omega^{-2} D\cap  (D+1)|+| D\cap  (\omega^{-2} D+1)|+4. 
\end{equation}

Let
\begin{align*}
&N_1=|\omega D\cap (\omega^3 D+1)|+|\omega^3 D\cap  (\omega D+1)|,\quad 
N_2=|\omega^{-1} D\cap (\omega D+1)|+|\omega D\cap  (\omega^{-1} D+1)|,\\
&N_3=|D\cap (\omega^2 D+1)|+|\omega^2 D\cap  (D+1)|,\quad
N_4=|\omega^{-2} D\cap  (D+1)|+| D\cap  (\omega^{-2} D+1)|.
\end{align*}
Then, \eqref{eq:redu1} and \eqref{eq:redu2} are rewritten as $N_1=N_2$ and $N_3=N_4+4$, respectively.
From the definition of $I$ and Table~\ref{Tab1} of Proposition~\ref{prop:cyclo8}, we have 
\begin{align*}
N_1=&\,(1,3)_8+(1,5)_8+(1,6)_8+(3,3)_8+(3,5)_8+(3,6)_8+(4,3)_8+
(4,5)_8+(4,6)_8\\
&\, +
(3,1)_8+(5,1)_8+(6,1)_8+(3,3)_8+(5,3)_8+(6,3)_8+(3,4)_8+
(5,4)_8+(6,4)_8\\
=&\,8n_2+n_3+4n_5+2n_6+2n_7+n_8, \\
N_2=&\,(7,1)_8+(7,3)_8+(7,4)_8+(1,1)_8+(1,3)_8+(1,4)_8+(2,1)_8+
(2,3)_8+(2,4)_8\\
&\, +
(1,7)_8+(3,7)_8+(4,7)_8+(1,1)_8+(3,1)_8+(4,1)_8+(1,2)_8+
(3,2)_8+(4,2)_8,\\
=&\, 
8n_2+n_3+4n_5+2n_6+2n_7+n_8,\\
N_3=&\,(0,2)_8+(0,4)_8+(0,5)_8+(2,2)_8+(2,4)_8+(2,5)_8+(3,2)_8+
(3,4)_8+(3,5)_8\\
&\, +
(2,0)_8+(4,0)_8+(5,0)_8+(2,2)_8+(4,2)_8+(5,2)_8+(2,3)_8+
(4,3)_8+(5,3)_8,\\
=&n_1+4n_2+2n_3+n_4+2n_5+n_6+3n_7+4n_8, \\
N_4=&\,(6,0)_8+(6,2)_8+(6,3)_8+(0,0)_8+(0,2)_8+(0,3)_8+(1,0)_8+
(1,2)_8+(1,3)_8\\
&\, +
(0,6)_8+(2,6)_8+(3,6)_8+(0,0)_8+(2,0)_8+(3,0)_8+(0,1)_8+
(2,1)_8+(3,1)_8\\
=&2n_1+6n_2+4n_3+2n_5+2n_6+2n_8. 
\end{align*}
It is clear that $N_1=N_2$. 
By the evaluations for $n_1,n_2,\ldots,n_8$ in Proposition~\ref{prop:cyclo8}, 
we have $N_3=(18q^2+28q-8a-32b-46)/64$ and 
$N_4=(18q^2+20q+8a+32b-46)/64$. Hence, $N_3=N_4+4$ if and only if 
$3q=a+4b+16$. This shows that \eqref{eq:cccc} holds if $3q=a+4b+16$.  
\qed
\vspace{0.3cm}

It is not difficult to see that the condition $q^2=a^2+2b^2$ with $3q=a\pm 4b+16$ and $a\equiv 1\,(\mod{4})$ is equivalent to that $q$ has the form $q=12c^2+4c+3$ with $c$ an arbitrary integer; in this case,  $a=4c^2+12c+1$ and $b=\pm (8c^2-2)$.  
Hence, by Theorem~\ref{thm:b2b3} and Proposition~\ref{prop:b0b2},  Theorem~\ref{thm:mainth} now follows. 

To see whether we have constructed an infinite family of Hadamard matrices in Theorem~\ref{thm:mainthm}, a natural question arises: are there infinitely many prime powers $q$ of the form $q=12c^2+4c+3$ with $c$ an integer? We believe that there are infinitely many primes of the form $12c^2+4c+3$ with $c$ an integer. But this is probably very difficult to prove. On the other hand, we conjecture that there are no proper prime powers $q$ of the form $q=12c^2+4c+3$ ($c$ is an integer). That is, we conjecture that there are no solutions to the equation
$$12c^2+4c+3=p^{\alpha},\; \alpha>1,$$
where $c$ is an integer, and $p$ is a prime. Some evidence is given in Introduction, namely all $58$ prime powers listed in \eqref{eq:list} are actually primes.  %(Question for Koji, did you really search for prime powers less than $10^5$ in Example 3.6, or you just searched for PRIMES less than $10^5$? I found it amazing that all 58 numbers listed in Example 3.6 are actually all primes; I believe that the Diophantine equation above should have only finitely many solutions. But the data in Example 3.6 suggests that it has no solutions at all).

\vspace{0.3cm}

%\section*{Concluding remarks}
%Let $G_0$ be any multiplicatively written  abelian  group of odd order, and define the group $G=\{x,xy:x\in G_0, y^2=1,yxy=x^{-1}\}$. Let  $A,B\subseteq G_0$ and $X=A \cup By$. Since $(xy)^{-1}=xy$ for $x\in G_0$, we have 
%\[
%XX^{(-1)}=AA^{(-1)}+BB^{(-1)}+(AB+BA)y. 
%\] 
%Hence, for  $A_i,B_i\subseteq G_0$, $i=0,1,2,3$, $\{\{0\}\times A_i\cup \{1\}\times B_i: i=0,1,2,3\}$ forms a difference family in $\Z_2\times G_0$ if and only if so does $\{A_i\cup B_i^{(-1)}y: i=0,1,2,3\}$ in $G$. 

%%%%%%%%%%%%%%%%%%%%%%%%%%%%%%%%%%%%%%%%%%%%%%%%%%%%%%%%%%%%
%%%%%%%%%%%%%%%%%%%%%%%%%%%%%%%%%%%%%%%%%%%%%%%%%%%%%%
%%%%%%%%%%%%%%%%%%%%%%%%%%%%%%%%%%%%%%%%%%%%%%%%%%%%%%%%%%


\begin{thebibliography}{99}

%\bibitem{BEW97}
%B. Berndt, R. Evans, K. S. Williams, {\it Gauss and Jacobi Sums}, Wiley, 1997. 

%\bibitem{bjl} T. Beth, D. Jungnickel, H. Lenz, {\it Design Theory}, Vol. I, 2nd edit., Cambridge University Press, 
%Cambridge, 
%1999.

%\bibitem{LN97}
%R. Lidl, H. Niederreiter, {\it Finite Fields}, Cambridge
%Univ. Press, 1997.

%\bibitem{BWX99} 
%A. E. Brouwer, R. M. Wilson, Q. Xiang, Cyclotomy and strongly regular graphs, {\it J. Alg. Combin.}, {\bf 10} (1999),
%25--28.

\bibitem{GS67}
J.-M. Goethals, J. J. Seidel, Orthogonal matrices with 
zero diagonal,  {\it Canad. J. Math.} 
{\bf 19} (1967), 1001--1010. 

\bibitem{K96}
H. Kimura, Hadamard matrices and dihedral groups, {\it Des. Codes Cryptogr.} {\bf 8} (1996), 71--77.  

\bibitem{KN02}
H. Kimura, T. Niwasaki, Some properties of Hadamard matrices coming from dihedral groups, {\it Graphs Combin.} {\bf 8} (2002), 319--327.  

\bibitem{LMS06}
K. H. Leung, S. L. Ma, B. Schmidt, 
New Hadamard matrices of order $4p^2$ obtained from Jacobi sums 
of order $16$,   
{\it J. Combin. Theory, Ser. A} {\bf 113} (2006), 822--838. 

\bibitem{LM18}
K. H. Leung, K. Momihara, New constructions of Hadamard matrices, {\tt arXiv:1809.05253}. 

\bibitem{peisert}
W. Peisert, All self-complementary symmetric graphs, {\it J. Algebra} {\bf 240} (2001), 209--229.


\bibitem{SY00}
K. Shinoda, M. Yamada, A family of Hadamard matrices of dihedral group type, 
 {\it Discrete Appl. Math.}  {\bf 102} (2000), 141--150. 

\bibitem{S75}
E. Spence, Hadamard matrices from relative difference sets, 
 {\it J. Combin. Theory, Ser. A} {\bf 19} (1975), 287--300. 


\bibitem{S67}
T. Storer, {\it Cyclotomy and Difference Sets},  Markham Publishing Company, 1967.

%\bibitem{T65}
%R. J. Turyn, Character sums and difference sets, {\it Pacific J. Math.,}
%{\bf 15} (1965), 319--346.


%\bibitem{Yamamoto}
%K. Yamamoto, On Jacobi sums and difference sets, {\it J. Combin. Theory, Ser. A,} {\bf 3} (1967), 146--181. 

\bibitem{WSW}
W. D. Wallis, A. P. Street, J. S. Wallis, {\it Combinatorics: Room Squares, 
Sum-Free Sets, Hadamard Matrices}, Lecture Notes in Mathematics, {\bf 292}, Springer, New York, 1972.

%\bibitem{W71}
%A. L. Whiteman, An infinite family of skew Hadamard matrices,  
%{\it Pacific J. Math.} {\bf 38} (1971), 817--822.

\bibitem{W72}
A. L. Whiteman, Hadamard matrices of order $4(2p+1)$,  
{\it Notices Amer. Math. Soc.} {\bf 19} (1972), A-681. 


%\bibitem{X92}
%M.-Y. Xia, Some infinite classes of special 
%Williamson matrices and difference sets,  {\it J. Combin. Theory, Ser. A,} {\bf 61} (1992), 230-242. 

\bibitem{XL91}
M.-Y. Xia, G. Liu, An infinite class of supplementary 
difference sets and Williamson matrices, {\it J. Combin. Theory, Ser. A} {\bf 58} (1991), 310--317. 

\bibitem{XL95}
M.-Y. Xia, G. Liu, On the class ${\mathcal H}_1^\ast$, 
{\it Acta Math. Sci.} {\bf 15} (1995), 361--369. 

\bibitem{XL96}
M.-Y. Xia, G. Liu, A new family of supplementary difference sets and Hadamard matrices, {\it J. Statist. Plann. Inference} {\bf 51} (2003), 263--275. 

\bibitem{XX94}
M.-Y. Xia, T. B. Xia, Hadamard matrices constructed from  supplementary 
difference sets in the class ${\mathcal H}_1$,  {\it J. Combin. Des.} 
{\bf 2} (1994), 325--339. 

\bibitem{XX99}
M.-Y. Xia, T. B. Xia, A family of $C$-partitions and $T$-matrices, 
 {\it J. Combin. Des.} 
{\bf 7} (1999), 269--281. 

\bibitem{XXSW05}
M.-Y. Xia, T. B. Xia, J. Seberry, J. Wu, An infinite 
family of Goethals-Seidel arrays,  {\it Discrete Appl. Math.}  
{\bf 145} (2005), 498--504. 

\bibitem{xiang}
Q. Xiang, Difference families from lines and half lines, {\it Europ. J. Combin.} {\bf 19} (1998), 395--400.

\end{thebibliography}
\end{document}